\documentclass[11pt]{article}
\usepackage{color}
\usepackage{amsfonts}
\newtheorem{theo}{Theorem}[section]
\newtheorem{lem}[theo]{Lemma}
\newtheorem{prop}[theo]{Proposition}

\newcommand{\mysection}[1]{\section{#1} \setcounter{equation}{0}}
\newcommand{\dfrac}[2]{\displaystyle \frac{#1}{#2}}

\newcommand{\vs}{\medskip}
\newcommand{\proof}{{\sc Proof.} \quad}

\newcommand{\R}{\mathbb{R}}

\newcommand{\be}{\begin{equation} \label}
\newcommand{\ee}{\end{equation}}
\newcommand{\bes}{\begin{equation} \begin{array}{c} \label}
\newcommand{\ees}{\end{array} \end{equation}}
\newcommand{\bea}{\begin{eqnarray}\label}
\newcommand{\eea}{\end{eqnarray}}
\newcommand{\bas}{\begin{eqnarray*}}
\newcommand{\eas}{\end{eqnarray*}}
\newcommand{\bit}{\begin{itemize}}
\newcommand{\eit}{\end{itemize}}
\newcommand{\nn}{\nonumber}
\newcommand{\eps}{\varepsilon}
\newcommand{\abs}{\\[3mm]}
\newcommand{\parab}{{\cal{P}}}

\newcommand{\ov}{\overline{v}}
\newcommand{\uv}{\underline{v}}

\def\qed{\hfill 
\raise -2pt
\hbox{\vrule \vbox to8pt{\hrule width 6pt
\vfill\hrule}\vrule}\\[5mm]}
\newcommand{\tg}{\vartheta}
\begin{document}
\title{\bf Rate of Convergence to Separable Solutions \\
of the Fast Diffusion Equation  }
\author{Marek Fila \\
Department of Applied Mathematics and Statistics, Comenius
University \\
84248 Bratislava, Slovakia \\
\\
and \\
\\
Michael Winkler \\
Institut f\"ur Mathematik, Universit\"at Paderborn \\
33098 Paderborn, Germany \\
}
\date{}
\maketitle
\begin{abstract}
\noindent
 We study the asymptotic behaviour near extinction of positive
solutions of the Cauchy problem
for the fast diffusion equation with a subcritical exponent.
We show that separable solutions are stable in some suitable sense
by finding a class of functions which belong to their domain of attraction.
For solutions in this class we establish optimal rates of convergence to 
separable solutions.

\end{abstract}
\mysection{Introduction}
We consider the Cauchy problem for the fast diffusion equation,
\begin{equation}\label{0}
       \left\{ \begin{array}{ll}
       u_\tau = \Delta (u^m),
               \qquad & x\in\R^n, \ \tau\in (0,T), \\[2mm]
       u(x,0)=u_0(x)\ge 0, \qquad & x\in\R^n,
       \end{array} \right.
\end{equation}
where $n\ge 3$, $T>0$ and $0<m<1$.
It is known that for $m<m_c:=(n-2)/n$ all
solutions with initial data satisfying
\[
u_0(x)=O\left(|x|^{-\frac{2}{1-m}}\right)\qquad\mbox{as}\quad |x|\to\infty,
\]
extinguish in finite time. We shall
consider solutions which vanish at $\tau=T$
and study their behaviour near $\tau=T$.

The function
\be{sep}
u(x,\tau):=((1-m)(T-\tau))^{\frac{1}{1-m}}\varphi^{\frac{1}{m}}(x)
\ee
is a solution of the fast diffusion equation $u_\tau = \Delta (u^m)$ if
$\varphi$ satisfies
\be{es}
    \Delta\varphi +\varphi^p=0, \qquad x\in\R^n, \qquad p:=\frac{1}{m} \, .
\ee
We call a nontrivial solution of the form (\ref{sep}) separable. We shall
show that separable solutions are stable in a suitable sense if
\[
n>10,\qquad 0<m<\frac{(n-2)(n-10)}{(n-2)^2-4n+8\sqrt{n-1}}\, .
\]
We also find optimal rates of convergence to separable solutions.

To formulate our results in more detail, it is convenient to introduce the
following change of variables:
\[
v(x,t):=((1-m)(T-\tau))^{-\frac{m}{1-m}}u^m(x,\tau),
\qquad t:=-\frac{1}{1-m}\ln\frac{T-\tau}{T}\, .
\]

If $u$ is a solution of (\ref{0}) with extinction at $\tau=T$ then $v$
satisfies
\begin{equation}  \label{eq:main}
 \left \{
\begin{array}{ll}
     (v^p)_t=\Delta v+v^p,\qquad & x\in \R^n,\  t>0,\qquad p=\frac{1}{m}\, ,
\vspace{5pt} \\
      v(x,0) =v_0(x), & x\in \R^n,
\end{array} \right.
\end{equation}
where
\[
v_0(x):=((1-m)T)^{-\frac{m}{1-m}}u_0^m(x).
\]
The behaviour of $u$ as $\tau\to T$ corresponds to the behaviour of $v$ as
$t\to\infty$. Hence, asymptotically separable extinction of $u$ is
equivalent to convergence of $v$ to a steady state as $t\to\infty$.

Concerning the existence of positive solutions of (\ref{es}),
it is well known that the Sobolev exponent
\[  p_S:=\left\{\
    \begin{array}{cl}
     \dfrac{{n+2}}{n-2}& \mbox{ if }
       n>2, \vs \\
       \infty & \mbox{ if }
       n\le 2,
    \end{array}\right.
\]
plays a crucial role.  Namely, there is a family of positive radial
solutions of (\ref{es})
if and only if $p\ge p_S$.  We denote the solution by
$\varphi=\varphi_\alpha(r)$, $r=|x|$, $\alpha>0$,  where
$\varphi_\alpha(r)$ satisfies
\begin{equation}  \label{eq:varphi}
\left \{
\begin{array}{l}
(\varphi_{\alpha})_{rr}+\dfrac{n-1}{r} (\varphi_{\alpha})_{r}+\varphi_\alpha^p=0,
\qquad r>0,\vs \\
\varphi_\alpha(0)=\alpha, \quad (\varphi_{\alpha})_{r}(0)=0.
\end{array} \right.
\end{equation}
For each $\alpha>0$, the solution $\varphi_\alpha$ is decreasing in
$|x|$ and
satisfies $\varphi_\alpha(|x|) \to 0$ as $|x| \to \infty$.

For the structure of the set of steady states,
there is another critical exponent
 \[
  p_c:=\left\{\
    \begin{array}{cl}
       \dfrac{(n-2)^2-4n+8\sqrt{n-1}}{(n-2)(n-10)}& \mbox{ if }
       n>10,  \vs \\
       \infty & \mbox{ if }
       n\le 10.
    \end{array}\right.
 \]
It is known that for $p_S\leq p<p_c$,
any positive steady state  intersects with other positive steady
states (see \cite{Wang}).  For $p\ge p_c$,
Wang \cite{Wang} showed that
the family of steady
states $\{\varphi_\alpha\, ; \, \alpha \in \R\}$
is ordered, that is, $\varphi_\alpha$ is 
increasing in $\alpha$ for each $x$. Moreover,
\[
  \lim_{\alpha\to 0}\varphi_\alpha(|x|)=0,\qquad \lim_{\alpha\to
    \infty}\varphi_\alpha (|x|)=\varphi_\infty(|x|),
\]
where $\varphi_\infty$ is a singular steady state given by
\[
\varphi_\infty(|x|):=L|x|^{-\nu}, \quad x \in \R^n\setminus\{0\},
\]
with
\be{L}
   \nu:=\dfrac{2}{p-1}\, ,   \qquad L:=\Big\{
\nu\left(n-2-\nu\right)\Big\}^{1/(p-1)}.
\ee
Each positive steady state satisfies
\begin{equation}  \label{phi_exp}
   \varphi_\alpha(|x|) = L|x|^{-\nu} - a_\alpha   |x|^{-\nu-\lambda_1} + 
o\Big(|x|^{-\nu-\lambda_1}\Big)
        \qquad \mbox{as } |x| \to\infty,
\end{equation}
where $\lambda_1$ is a positive constant given by
\[
\lambda_1 =\lambda_1 (n,p)
  := \frac{n-2-2\nu-\sqrt{(n-2-2\nu)^2-8(n-2-\nu)}}{2}\,,
\]
and  $a_\alpha>0$ is such that
\be{a_alpha}
a_\alpha>a_\beta\qquad \mbox{if }0<\alpha<\beta.
\ee 
(see  \cite{Gui-N-W:stab}).

%
%
%
%
%
Our main results are contained in the following two theorems:
\begin{theo}\label{theo100}
  Let $n>10$, $p>p_c$, $\alpha>0$, and assume that $v_0$ is continuous in $\R^n$ and satisfies
\be{sing}
        0 \le v_0(x) \le L|x|^{-\nu},\qquad x \in \R^n\setminus\{0\},
\ee
and
  \be{100.1}
	|v_0(x)-\varphi_\alpha(|x|)| \le b |x|^{-\gamma}, \qquad
x\in\R^n,\quad |x|>1,
  \ee
  with some $b>0$ and 
  \be{gamma}
	\gamma \in \Big(\nu+\lambda_1 \, , \, \frac{n-2}{2} \Big).
  \ee
  Then there exists $C>0$ such that the solution $v$ of {\rm (\ref{eq:main})} has the property that
  \be{100.3}
	|v(x,t)-\varphi_\alpha(|x|)| \le C e^{-\kappa(\gamma) t},
	\qquad x\in\R^n, \quad t\ge 0,
  \ee
  where the positive number $\kappa(\gamma)$ is given by 
  \be{kappa_gamma}
	\kappa(\gamma):=\frac{\gamma(n-2-\gamma)}{pL^{p-1}}-1
	=\frac{\gamma(n-2-\gamma)}{(\nu+2)(n-2-\nu)}-1\, .
  \ee
\end{theo}
The convergence rate provided by Theorem~\ref{theo100} is in fact optimal for stabilization from above
or from below:
\begin{theo}\label{theo200}
  Let $n>10$, $p>p_c$, $\alpha>0$ and $b>0$, and assume that $\gamma$ and $\kappa(\gamma)$ are as in
{\rm (\ref{gamma})} and {\rm (\ref{kappa_gamma})},
  respectively.
  Moreover, let $v$ denote the solution of {\rm (\ref{eq:main})} corresponding to some
continuous $v_0$. \abs
%
  {\rm (i)} \ If
  \be{200.1}
	v_0(x) \ge \varphi_\alpha(|x|)+b(|x|+1)^{-\gamma},
	\qquad x\in\R^n,
  \ee
  then there exists $C>0$ such that
  \be{200.2}
	v(0,t) - \varphi_\alpha(0) \ge C e^{-\kappa(\gamma)t},
	\qquad t\ge 0.
  \ee
  {\rm (ii)} \ Under the assumption that
  \be{200.3}
	v_0(x) \le \varphi_\alpha(|x|)-b(|x|+1)^{-\gamma},
	\qquad x\in\R^n,
  \ee
  one can find $C>0$ such that
  \be{200.4}
	\varphi_\alpha(0) - v(0,t) \ge C e^{-\kappa(\gamma)t},
	\qquad t\ge 0.
  \ee
\end{theo}

Theorems~\ref{theo100} and \ref{theo200} imply that there is a continuum of
rates of convergence of solutions of (\ref{eq:main}) to $\varphi_\alpha$,
$\alpha>0$, since the range of $\kappa$ is
\[
\left(0,\frac{(n-2)^2}{4pL^{p-1}}-1\right).
\]
Notice that $\kappa(\gamma)\to 0$ as $\gamma\to\nu+\lambda_1$.

Asymptotically separable extinction of solutions of (\ref{0}) has only been
known to occur for $m=(n-2)/(n+2)=1/p_S$, $n>2$. This was predicted in \cite{K1} and then
rigorously established in \cite{GP}, \cite{dPS}, Theorem~7.10 in \cite{Vsmooth}. 
The exponent $m=1/p_S$ is the unique value of $m$ for which
the second-kind selfsimilar solution is separable, see \cite{K1}, \cite{K2},
\cite{PZ}. The rate of convergence to a separable solution is not known in
this case.

Depending on the decay rate of the initial function $u_0$, other types of
asymptotic behaviour near extinction may occur, such as convergence to
Barenblatt profiles (see \cite{BBDGV}, \cite{BDGV}, \cite{BGV}, \cite{FKW},
\cite{FVW}, \cite{FVWY}, for example) or convergence to selfsimilar
solutions of the second type (see \cite{GP}, \cite{Vsmooth}).

For a smoothly bounded domain $\Omega\subset\R^n$ all bounded positive
solutions of $u_\tau = \Delta (u^m)$ with the homogeneous Dirichlet boundary
condition and $(n-2)_+/(n+2)<m<1$ extinguish in finite time and they
approach separable solutions, see \cite{AK}, \cite{BH}, \cite{BGVb},
\cite{Kw}, \cite{SV}. But little is known about the convergence rate. As far
as we know, only upper bounds for the decay rates of the entropy and of a
weighted $L^2$-norm of the relative distance from the separable solution are
given in \cite{BGVb} for $m$ sufficiently close to 1.

Let us also compare Theorems~\ref{theo100} and \ref{theo200} with corresponding
results from \cite{FWY}, \cite{HY} on convergence to steady states for the
problem
\begin{equation}  \label{eq:fujita}
 \left \{
\begin{array}{ll}
     v_t=\Delta v+v^p,\qquad & x\in \R^n,\  t>0,\qquad p>p_c,
\vspace{5pt} \\
      v(x,0) =v_0(x), & x\in \R^n.
\end{array} \right.
\end{equation}
For both problems, the steady states are the same and we obtain a continuum
of convergence rates which depend explicitly on the tail of initial
functions. But for (\ref{eq:main}) the rates are exponential while for
(\ref{eq:fujita}) they are algebraic, see \cite{FWY}, \cite{HY}.

Similarly as for (\ref{eq:fujita}), the steady states $\varphi_\alpha$
are unstable from below and from above for (\ref{eq:main})
when $p_S\le p<p_c$.

\begin{prop}\label{prop:unst}
  Let $n>2$, $p_S\le p<p_c$, $\alpha>0$, and assume that $v_0$ is bounded
and continuous in
$\R^n$.\abs
  {\rm (i)} \ If
\[
        v_0(x) \ge \varphi_\alpha(|x|), \qquad x\in\R^n,
\qquad v_0\not\equiv\varphi_\alpha \, ,
\]
  then
\[
\Vert v(\cdot,t)\Vert_{L^\infty(\R^n)}
\to\infty \qquad\mbox{as}\quad t\to \infty.
\]
  {\rm (ii)} \ If
\[
        0\le v_0(x) \le \varphi_\alpha(|x|), \qquad x\in\R^n,
\qquad v_0\not\equiv\varphi_\alpha \, ,
\]
then there is $T\in(0,\infty]$ such that
\[
\Vert v(\cdot,t)\Vert_{L^\infty(\R^n)}\to 0 \qquad\mbox{as}\quad t\to T.
\]
\end{prop}

The paper is organised as follows. In Section~2 we collect some
preliminaries.
In Section~3 we prove Theorem~\ref{theo100} for radially
symmetric solutions with $v_0\ge \varphi_\alpha$. In Section~4 we give a 
corresponding bound for $v_0\le \varphi_\alpha$ and we complete the proof of
Theorem~\ref{theo100}. Sections~5 and 6 are devoted to the proof of
Theorem~\ref{theo200} and Section~7 to Proposition~\ref{prop:unst}.
\mysection{Preliminaries}
\subsection{An ODE lemma}
The following statement plays a key role in our analysis.
It describes the spatial profile of separated solutions to the formal linearization of 
(\ref{eq:main}) around
a given steady state $\varphi_\alpha$.
\begin{lem}\label{lem1}
  Let $\alpha>0$ and $\kappa$ be positive and such that
  \be{1.1}
	\kappa<\kappa_0(p):=\frac{(n-2)^2}{4pL^{p-1}}-1.
  \ee
  Then the solution $f=f_{\alpha,\kappa}$ of the initial-value problem
  \be{1.2}
	\left\{ \begin{array}{l}
	f_{rr} + \frac{n-1}{r} f_r + (\kappa+1) p \varphi_\alpha^{p-1} f=0, \qquad r>0, \\[1mm]
	f(0)=1, \quad f_r(0)=0,
	\end{array} \right.
  \ee
  is positive on $[0,\infty)$, and there exist $c>0$ and $C>0$ such that
  \be{1.3}
	cr^{-\gamma} \le f(r) \le Cr^{-\gamma}
	\qquad \mbox{for all } r>1,
  \ee
  where $\gamma=\gamma(\kappa)$ is the positive number given by
  \be{1.4}
	\gamma(\kappa):=\frac{n-2}{2} - \sqrt{\frac{(n-2)^2}{4} - (\kappa+1)pL^{p-1}}.
  \ee
\end{lem}
\proof
  We first note that as an evident consequence of (\ref{1.1}) and the fact that $n\ge 3$, the 
number $\gamma=\gamma(\kappa)$
  given by (\ref{1.4}) indeed is real and positive.
  Since $f(0)=1$, we know that
  \bas
	r_\star:=\sup \Big\{ r>0 \ \Big| \ f>0 \mbox{ on } [0,r] \Big\}
  \eas
  is a well-defined element of $(0,\infty]$, and for $\tg>0$ we let
  \bas
	h(r) \equiv h_{\tg}(r) := r^{\tg} f(r), \qquad r\in [0,r_\star).
  \eas
  Then 
  \be{1.44}
	h(0)=0
	\qquad \mbox{and} \qquad
	h>0 \mbox{ on } (0,r_\star),
  \ee
  and since 
  \bas
	& & f_r(r)=r^{-\tg} h_r(r) - \tg r^{-\tg-1} h(r), \quad r\in (0,r_\star), 
	\qquad \mbox{and} \\
	& & f_{rr}(r)=r^{-\tg} h_{rr}(r) - 2\tg r^{-\tg-1} h_r(r) + \tg(\tg+1) r^{-\tg-2} h(r), \quad r\in (0,r_\star),
  \eas
  from (\ref{1.2}) we obtain that
  \bea{1.5}
	0 &=& r^{-\tg} h_{rr} - 2\tg r^{-\tg-1} h_r + \tg(\tg+1) r^{-\tg-2} h \nn\\
	& & + \frac{n-1}{r} \Big\{ r^{-\tg} h_r - \tg r^{-\tg-1} h \Big\} 
	 + (\kappa+1) p \varphi_\alpha^{p-1} 
 r^{-\tg} h \nn\\[1mm]
	&=& r^{-\tg} \Big\{ h_{rr} + \frac{n-1-2\tg}{r} h_r 
	+ \frac{\tg(\tg+2-n)+(\kappa+1)p\varphi_\alpha^{p-1} r^2}{r^2} h \Big\} \nn\\[1mm]
  \eea
for all $r\in (0,r_\star)$.
  We first apply this to $\tg:=\gamma$ and $h:=h_\gamma$ and use that 
\[\varphi_\alpha(r) \le Lr^{-\nu}\]
  for all $r>0$, which by positivity of $h$ on $(0,r_\star)$ allows us to conclude that
  \be{1.6}
	0 \le h_{rr} + \frac{n-1-2\gamma}{r} h_r + \frac{\gamma(\gamma+2-n)+(\kappa+1)pL^{p-1}}{r^2} h
	\qquad \mbox{for all } r\in (0,r_\star).
  \ee
  But the definition (\ref{1.4}) entails that $\gamma$ coincides with the smaller root of the equation
  \be{1.66}
	\gamma(\gamma+2-n) + (\kappa+1)pL^{p-1}=0,
  \ee
  whence (\ref{1.6}) actually says that
  \bea{1.7}
	0 \le h_{rr} + \frac{n-1-2\gamma}{r} h_r 
	= r^{-n+1+2\gamma} \Big(r^{n-1-2\gamma} h_r \Big)_r
\eea
for all $r\in (0,r_\star)$.
  Now according to (\ref{1.44}) it is clearly possible to fix $r_1\in (0,r_\star)$ small enough such that $r_1<1$
  and $h_r(r_1) \ge 0$. Therefore, integrating (\ref{1.7}) over $(r_1,r)$ for arbitrary $r\in (r_1,r_\star)$ yields
  \bas
	r^{n-1-2\gamma} h_r(r) \ge r_1^{n-1-2\gamma} h_r(r_1) \ge 0
	\qquad \mbox{for all } r\in (r_1,r_\star),
  \eas
  so that
  \bas
	h(r) \ge c_1
	\qquad \mbox{for all } r\in (r_1,r_\star),
  \eas
  where $c_1:=h(r_1)$ is positive thanks to (\ref{1.44}). 
  In particular, this ensures that necessarily $r_\star=\infty$, and moreover we obtain that
  \bas
	f(r)=r^{-\gamma}h(r) \ge c_1 r^{-\gamma}
	\qquad \mbox{for all } r>r_1,
  \eas
  which implies the left inequality in (\ref{1.3}), because $r_1<1$.

  Our derivation of the upper estimate for $f$ in (\ref{1.3}) will proceed in two steps.
  First, since $\gamma<(n-2)/2$ by (\ref{1.4}) and (\ref{1.1}), we can pick $\delta>0$ such that
  \be{1.777}
	\delta<n-2-2\gamma
  \ee
  and $\delta \le \lambda_1$. We then fix $\eps>0$ small such that
  \be{1.77}
	\gamma_\eps:=\frac{n-2}{2} - \sqrt{\frac{(n-2)^2}{4} - (\kappa+1)(pL^{p-1}-\eps)}
  \ee
  satisfies
  \be{1.8}
	\gamma_\eps>\gamma-\delta,
  \ee
  which can clearly be achieved by a continuity argument.
  Now as a consequence of (\ref{phi_exp}) we can find $c_2>0$ such that
  \be{1.9}
	\varphi_\alpha(r) \ge Lr^{-\nu} - c_2 r^{-\nu-\lambda_1}
	\qquad \mbox{for all } r>0,
  \ee
  whence in particular $\varphi_\alpha^{p-1}(r) r^2 \to L^{p-1}$ as $r\to\infty$. Accordingly,
  there exists $r_2>0$ such that
  \bas
	p\varphi_\alpha^{p-1}(r) r^2 \ge pL^{p-1}-\eps
	\qquad \mbox{for all } r \ge r_2,
  \eas
  so that (\ref{1.5}) applied to $\tg:=\gamma_\eps$ shows that $h:=h_{\gamma_\eps}$ satisfies
  \bas
	0 \ge h_{rr} + \frac{n-1-2\gamma_\eps}{r} h_r
	+ \frac{\gamma_\eps(\gamma_\eps+2-n)+(\kappa+1)(pL^{p-1}-\eps)}{r^2} h
	\qquad \mbox{for all } r\ge r_2.
  \eas
  Due to (\ref{1.77}), the zero-order term again vanishes, and hence we have
  \bas
	0 \ge h_{rr} + \frac{n-1-2\gamma_\eps}{r^2} h_r
	= r^{-n+1+2\gamma_\eps} (r^{n-1-2\gamma_\eps} h_r)_r
	\qquad \mbox{for all } r\ge r_2,
  \eas
  on successive integration implying that
  \bas
	r^{n-1-2\gamma_\eps} h_r(r) \le c_3:= r_2^{n-1-2\gamma_\eps} h_r(r_2)
	\qquad \mbox{for all } r\ge r_2
  \eas
  and
  \bas
	h(r) \le h(r_2) + c_3 \int_{r_2}^r \rho^{1-n+2\gamma_\eps} d\rho
	\qquad \mbox{for all } r\ge r_2.
  \eas
  Since from (\ref{1.77}) we know that $\gamma_\eps < (n-2)/2$ and thus $1-n+2\gamma_\eps<-1$, from this and 
  the definition of $h=h_{\gamma_\eps}$ we infer that
  \be{1.10}
	f(r)=r^{-\gamma_\eps} h(r) \le c_4 r^{-\gamma_\eps}
	\qquad \mbox{for all } r\ge r_2
  \ee
  with 
\[c_4:=h(r_2)+ \frac{c_3 r_2^{2-n+2\gamma_\eps}}{n-2-2\gamma_\eps}.\]

  In order to finally derive the second inequality in (\ref{1.3}) from this, we fix $c_5>0$ and $r_3>r_2$ large enough
  fulfilling
  \bas
	(1-z)^{p-1} \ge 1-c_5 z
	\qquad \mbox{for all } z\in [0,1/2]
\qquad \mbox{and}\qquad
	\frac{c_2}{L} r_3^{-\lambda_1} \le \frac{1}{2}.
  \eas
  Then (\ref{1.9}) says that for all $r\ge r_3$ we can estimate
  \bas
	p\varphi_\alpha^{p-1}(r) r^2
	&\ge& pr^2 \Big(Lr^{-\nu} - c_2 r^{-\nu-\lambda_1}\Big)^{p-1} 
	= pL^{p-1} \Big(1-\frac{c_2}{L} r^{-\lambda_1}\Big)^{p-1} \\
	&\ge& pL^{p-1} \Big(1-c_5 \frac{c_2}{L} r^{-\lambda_1}\Big) 
	= pL^{p-1} - c_6 r^{-\lambda_1}
  \eas
  with an evident definition of $c_6$.
  Hence, returning to our original choice $h=h_\gamma$, from (\ref{1.5}) and (\ref{1.66}) we obtain that
  \bea{1.11}
	0 &\ge& h_{rr} + \frac{n-1-2\gamma}{r} h_r 
	+ \frac{\gamma(\gamma+2-n)+(\kappa+1) (pL^{p-1}-c_6 r^{-\lambda_1})}{r^2} h \nn\\
	&=& r^{-n+1+2\gamma} (r^{n-1-2\gamma} h_r)_r - c_6 r^{-2-\lambda_1}h
	\qquad \mbox{for all } r\ge r_3.
  \eea
  Here we use (\ref{1.10}) and the fact that $\delta \le \lambda_1$ in estimating
  \bas
	c_6 r^{-2-\lambda_1} h(r)
	= c_6 r^{-2-\lambda_1+\gamma}f(r) 
	\le c_4 c_6 r^{-2-\lambda_1+\gamma-\gamma_\eps} 
	\le c_7 r^{-2-\delta+\gamma-\gamma_\eps}
	\qquad \mbox{for all } r\ge r_3
  \eas
  with $c_7:=c_4 c_6 r_3^{-(\lambda_1-\delta)}$, so that an integration in (\ref{1.11}) shows that if we abbreviate
  $c_8:=r_3^{n-1-2\gamma} h_r(r_3)$, then
  \be{1.12}
	r^{n-1-2\gamma} h_r(r)
	\le c_8 + c_7 \int_{r_3}^r \rho^{n-3-\gamma-\gamma_\eps-\delta} d\rho
	\qquad \mbox{for all } r\ge r_3.
  \ee
  Since the inequality $\gamma_\eps<\gamma$ along with (\ref{1.777}) guarantees that
  \bas
	n-2-\gamma-\gamma_\eps-\delta > n-3-2\gamma-\delta>0,
  \eas
  we have
  \bas
	\int_{r_3}^r \rho^{n-3-\gamma-\gamma_\eps-\delta}d\rho
	\le  \frac{r^{n-2-\gamma-\gamma_\eps-\delta}}{n-2-\gamma-\gamma_\eps-\delta}
	\qquad \mbox{for all } r\ge r_3,
  \eas
  so that (\ref{1.12}) implies that
  \bas
	h_r(r)
	\le c_8 r^{-n+1+2\gamma} + c_9 r^{-1+\gamma-\gamma_\eps-\delta}
	\qquad \mbox{for all } r\ge r_3, \qquad c_9:=\frac{c_7}{n-2-\gamma-\gamma_\eps-\delta}.
  \eas
 In view of our
restriction (\ref{1.8}) on $\gamma_\eps$,
  from this we obtain
  \bas
	h(r) \le h(r_3) + \frac{c_8}{n-2-2\gamma} + \frac{c_9}{\gamma_\eps-\gamma+\delta}
	\qquad \mbox{for all } r\ge r_3.
  \eas
  By definition of $h=h_\gamma$ and the positivity of $f$ on $[0,r_3]$, this finally completes the proof of (\ref{1.3}).
\qed
\subsection{Two auxiliary asymptotic estimates}
The following two auxiliary lemmata will be used both in Lemma~\ref{lem4} and in Lemma~\ref{lem8} in order to
provide appropriate control of certain higher order expressions arising during
linearization.
\begin{lem}\label{lem7}
  Let $\alpha>0$, $\beta>\alpha$, $\mu>0$ and $\kappa \in (0,\kappa_0(p))$, where $\kappa_0(p)$ is as defined in 
{\rm (\ref{1.1})}.
  Then for all $A\ge 0$ and $B\ge 0$ there exist $r_0>0$ and $c>0$ such that with $f_{\mu,\kappa}$ as in {\rm (\ref{1.2})}
  we have
  \be{7.1}
	\Big\{\varphi_\beta(r) - Bf_{\mu,\kappa}(r)\Big\}^{p-1} - \Big\{\varphi_\alpha(r) + Af_{\mu,\kappa}(r)\Big\}^{p-1} 
	\ge c r^{-2-\lambda_1}
	\qquad \mbox{for all } r>r_0.
  \ee
\end{lem}
\proof
  Since $\beta>\alpha$ and hence $a_\beta<a_\alpha$ by (\ref{a_alpha}), we can fix some small $\eta>0$ such that
  \be{7.2}
	c_1:=(p-1-\eta)(a_\alpha-\eta) - (p-1+\eta)(a_\beta+\eta) \, >0.
  \ee
  Then in view of (\ref{phi_exp}) we can find $r_1>1$ fulfilling
  \be{7.3}
	\varphi_\alpha(r) \le Lr^{-\nu} - \Big(a_\alpha-\frac{\eta}{2}\Big) r^{-\nu-\lambda_1}
	\qquad \mbox{for all } r>r_1
  \ee
  and
  \be{7.4}
	\varphi_\beta(r) \ge Lr^{-\nu} - \Big(a_\beta+\frac{\eta}{2}\Big) r^{-\nu-\lambda_1}
	\qquad \mbox{for all } r>r_1.
  \ee
  Now from Lemma~\ref{lem1} we know that there exists $c_2>0$ such that $f=f_{\mu,\kappa}$ satisfies
  \bas
	f(r) \le c_2 r^{-\gamma}
	\qquad \mbox{for all } r>1
  \eas
  with $\gamma=\gamma(\kappa)$ given by (\ref{1.4}). As $\gamma>\nu+\lambda_1$, we can therefore choose
  $r_2>r_1$ in such a way that
  \bas
	\max \{A,B\} f(r) \le \frac{\eta}{2} r^{-\nu-\lambda_1}
	\qquad \mbox{for all } r>r_2.
  \eas
  Together with (\ref{7.3}) and (\ref{7.4}), this yields the inequalities
  \be{7.5}
	\varphi_\alpha(r) + Af(r) \le Lr^{-\nu} - (a_\alpha-\eta) r^{-\nu-\lambda_1}
	\qquad \mbox{for all } r>r_2
  \ee
  and
  \be{7.6}
	\varphi_\beta(r) + Bf(r) \ge Lr^{-\nu} - (a_\beta+\eta) r^{-\nu-\lambda_1}
	\qquad \mbox{for all } r>r_2.
  \ee
  We next take $z_1>0$ small enough fulfilling
  \be{7.7}
	1-(p-1+\eta)z \le (1-z)^p \le 1-(p-1-\eta)z
	\qquad \mbox{for all } z\in [0,z_1],
  \ee
  and then fix $r_3>r_2$ such that
  \bas
	\frac{a_\alpha-\eta}{L} r_3^{-\lambda_1} \le z_1,
  \eas
  which by (\ref{7.2}) implies that also 
\[\frac{a_\beta+\eta}{L} r_3^{-\lambda_1} \le z_1.\] 
  Hence, (\ref{7.5}) and the second inequality in (\ref{7.7}) show that
  \bea{7.8}
	\Big\{\varphi_\alpha(r)+Af(r)\Big\}^{p-1}
	&\le& \Big(Lr^{-\nu}\Big)^{p-1} \Big\{1-\frac{a_\alpha-\eta}{L} r^{-\lambda_1}\Big\}^{p-1}\nn\\
	&\le& \Big(Lr^{-\nu}\Big)^{p-1} \Big\{1-(p-1-\eta) \frac{a_\alpha-\eta}{L} r^{-\lambda_1}
		\Big\} \nn\\
	&=& L^{p-1} r^{-2} - (p-1-\eta)(a_\alpha-\eta) L^{p-1} r^{-2-\lambda_1}
  \eea
  for all $r>r_3$, whereas similarly (\ref{7.6}) and the first inequality in (\ref{7.7}) entail that for any such $r$ we have
  \bea{7.9}
	\Big\{\varphi_\beta(r)-Bf(r)\Big\}^{p-1}
	&\ge& \Big(Lr^{-\nu}\Big)^{p-1} \Big\{1-(p-1+\eta) \frac{a_\beta+\eta}{L} r^{-\lambda_1}
		\Big\} \nn\\
	&=& L^{p-1} r^{-2} - (p-1+\eta)(a_\beta+\eta) L^{p-1} r^{-2-\lambda_1}.
  \eea
  In light of (\ref{7.2}), the desired estimate immediately results from (\ref{7.8}) and (\ref{7.9}).
\qed
\begin{lem}\label{lem9}
  Let $\alpha>0, \mu>0$ and $\kappa\in (0,\kappa_0(p))$ with $\kappa_0(p)$  from
{\rm (\ref{1.1})}. 
  Then there exists $C>0$ such that
  \be{9.1}
	\varphi_\alpha^{p-2}(r) f_{\mu,\kappa}(r) \le
C(r+1)^{\nu-\gamma-2}
	\qquad \mbox{for all } r\ge 0,
  \ee
  where $\gamma=\gamma(\kappa)$ is as in {\rm (\ref{1.4})}.
\end{lem}
\proof
  We apply (\ref{phi_exp}) and (\ref{1.3}) to find positive constants $c_1$ and $c_2$ such that
  \be{9.2}
	c_1 r^{-\nu} \le \varphi_\alpha(r) \le Lr^{-\nu}
	\qquad \mbox{for all } r>1,
  \ee
  and such that $f:=f_{\mu,\kappa}$ satisfies
  \be{9.3}
	f(r)\le c_2 r^{-\gamma}
	\qquad \mbox{for all } r>1.
  \ee
  Writing $c_3:=c_2 \max\{c_1^{p-2}, L^{p-2}\}$, we therefore see that
  \bas
	\varphi_\alpha^{p-2}(r) f(r) \le 
c_3 r^{\nu-\gamma-2}
	\le c_3 2^{\nu-\gamma-2} (r+1)^{\nu-\gamma-2}
	\qquad \mbox{for all } r>1,
  \eas
  because $\gamma>\nu+\lambda_1$ implies that
$\nu-\gamma-2<0$. 
  Since for $r\in [0,1]$, (\ref{9.1}) is obvious from the positivity of $\varphi_\alpha$ and the boundedness of both
  $\varphi_\alpha$ and $f$, the proof is complete.
\qed
\mysection{Convergence from above: upper bound}
In our first estimate of solutions from above, in Lemma \ref{lem3} below, we shall apply a comparison argument involving 
comparison functions which monotonically decrease with time. 
The initial data of the latter will be constructed separately in the following lemma.

\begin{lem}\label{lem2}
  Assume that $\alpha>0$ and $\kappa\in (0,\kappa_0(p))$ with $\kappa_0(p)$ as in {\rm (\ref{1.1})}, 
  and let $f_{\alpha,\kappa}$ denote the corresponding solution of {\rm (\ref{1.2})}.
  Then there exists $A_0=A_0(\alpha,\kappa)>0$ with the property that given any $A>A_0$ one can find $r_A>1$ such that
  for each $\eps\in (0,1)$, the number
  \be{2.1}
	r_{A\eps} := \sup \Big\{ \tilde r>0 \ \Big| \
	\varphi_\alpha(r)+Af_{\alpha,\kappa}(r) > L(r+\eps)^{-\nu} \quad \mbox{for all } r\in (0,\tilde r) \Big\}
  \ee
  is well-defined with
  \be{2.01}
	1 \le r_{A\eps} \le r_A,
  \ee
  and
  \be{2.2}
	v_{0A\eps}(r):=
	\left\{ \begin{array}{ll}
	L(r+\eps)^{-\nu} \qquad & \mbox{if } r\in (0,r_{A\eps}), \\[2mm]
	\varphi_\alpha(r)+Af_{\alpha,\kappa}(r) \qquad & \mbox{if } r\ge r_{A\eps},
	\end{array} \right.
  \ee
  determines a positive function 
  $v_{0A\eps} \in W^{1,\infty}_{loc}([0,\infty)) \cap C^2((0,\infty) \setminus \{r_{A\eps}\})$
  which satisfies
  \be{2.3}
	(v_{0A\eps})_{rr} + \frac{n-1}{r} (v_{0A\eps})_r + v_{0A\eps}^p \le 0
	\qquad \mbox{for all } r\in (0,\infty)\setminus \{r_{A\eps}\}
  \ee
  as well as
  \be{2.4}
	\liminf_{r\nearrow r_{A\eps}} (v_{0A\eps})_r(r)
	\ge
	\limsup_{r\searrow r_{A\eps}} (v_{0A\eps})_r(r).
  \ee
\end{lem}
\proof
  According to (\ref{phi_exp}), we can fix $c_1>0$ such that
  \be{2.44}
	\varphi_\alpha(r) \ge Lr^{-\nu} - c_1 r^{-\nu-\lambda_1}
	\qquad \mbox{for all }r>0.
  \ee
  Since $\kappa>0$, we can then take $z_1 \in (0,1)$ such that
  \be{2.5}
	(1+z)^p \le 1+(\kappa+1)pz
	\qquad \mbox{for all } z\in [0,z_1],
  \ee
  and pick $r_1>1$ large fulfilling
  \be{2.6}
	\frac{2c_1}{L} r_1^{-\lambda_1} \le z_1.
  \ee
  We thereupon choose $A_0>0$ large enough such that with $f:=f_{\alpha,\kappa}$ we have
  \be{2.7}
	\varphi_\alpha(r) + Af(r_1) > L r_1^{-\nu}
	\qquad \mbox{for all } A>A_0,
  \ee
  and let $A>A_0$ and $\eps\in (0,1)$ be given. 
  To see that then the set on the right-hand side of (\ref{2.1}) is nonempty and bounded, we first note that by (\ref{2.7})
  we have
  \bas
	\varphi_\alpha(r) + Af(r_1) > L (r_1+\eps)^{-\nu}.
  \eas
  Moreover, again by (\ref{phi_exp}) and by (\ref{1.3}) there exist $c_2>0$ and $c_3>0$ such that
  \bas
	\varphi_\alpha(r) \le Lr^{-\nu} - c_2 r^{-\nu-\lambda_1}
	\qquad \mbox{for all } r>1
  \eas
  and
  \bas
	f(r) \le c_3 r^{-\gamma(\kappa)}
	\qquad \mbox{for all } r>1
  \eas
  with $\gamma(\kappa)$ as in (\ref{1.4}). 
  Since the positivity of $\kappa$ implies that
  \bas
	\gamma(\kappa) > \frac{n-2}{2} - \sqrt{\frac{(n-2)^2}{4}-pL^{p-1}}
	= \nu+\lambda_1,
  \eas
  for some large $r_2(A)$ we thus obtain
  \be{2.8}
	\varphi_\alpha(r) + Af(r) \le Lr^{-\nu} - \frac{c_2}{2} r^{-\nu-\lambda_1}
	\qquad \mbox{for all } r>r_2(A).
  \ee
  On the other hand, by convexity of $[0,\infty) \ni z \mapsto (1+z)^{-\nu}$, we know that
  \be{2.9}
	L(r+\eps)^{-\nu}
	= Lr^{-\nu} \Big(1+\frac{\eps}{r}\Big)^{-\nu}
	\ge Lr^{-\nu} - \nu\eps L r^{-\nu-1}
	\qquad \mbox{for all } r>0.
  \ee
  Since it can easily be checked that $\lambda_1>1$, combining (\ref{2.8}) with (\ref{2.9}) we obtain $r_3(A)>0$
such that 
 \bas
	\varphi_\alpha(r) + Af(r) < L(r+\eps)^{-\nu}
	\qquad \mbox{for all } r>r_3(A).
  \eas
  Having thereby shown that $r_{A\eps}$ is well-defined and satisfies
  \be{2.10}
	1<r_1<r_{A\eps}<r_3(A),
  \ee
  we let $v_{0A\eps}$ be given by (\ref{2.2}) and proceed to verify (\ref{2.3}).
  To this end, for small $r$ we use the definition (\ref{L}) of $L$ in estimating
  \bas
	(v_{0A\eps})_{rr} + \frac{n-1}{r} (v_{0A\eps})_r + v_{0A\eps}^p
	= - \nu (n-1)L (r+\eps)^{-\frac{p+1}{p-1}} \Big(\frac{1}{r}-\frac{1}{r+\eps}\Big) 
	< 0
\eas
for all $ r<r_{A\eps}$.

  In the corresponding outer region, we first observe that from the definition of $r_{A\eps}$ it follows that
  \bas
	\varphi_\alpha(r) + Af(r) \le L(r+\eps)^{-\nu} 
	\qquad \mbox{for all }r>r_{A\eps},
  \eas
  so that according to (\ref{2.44}),
  \bea{2.11}
	Af(r) \le L(r+\eps)^{-\nu}
	- \left(Lr^{-\nu} - c_1 r^{-\nu-\lambda_1}\right) 
	\le c_1 r^{-\nu-\lambda_1}
	\qquad \mbox{for all }r>r_{A\eps}.
  \eea
  In conjunction with (\ref{2.6}) and (\ref{2.10}), (\ref{2.44}) furthermore guarantees that
  \bas
	\varphi_\alpha(r)
	\ge Lr^{-\nu} \Big(1-\frac{c_1}{L} r^{-\lambda_1}\Big) 
	\ge Lr^{-\nu} \Big(1-\frac{c_1}{L} r_1^{-\lambda_1}\Big) 
	\ge Lr^{-\nu} \Big(1-\frac{z_1}{2}\Big) 
	> \frac{L}{2} r^{-\nu}
\eas
for all $r>r_{A\eps}$,
  because $z_1<1$. Therefore, (\ref{2.11}) and again (\ref{2.6}) and (\ref{2.44}) yield
  \bas
	\frac{Af(r)}{\varphi_\alpha(r)} < \frac{2c_1}{L} r^{-\lambda_1} \le z_1
	\qquad \mbox{for all } r>r_{A\eps},
  \eas
  and hence (\ref{2.5}) applies to show that the nonlinear term in (\ref{2.3}) can be estimated according to
  \bas
	v_{0A\eps}^p(r) &=& \varphi_\alpha^p(r) \left(1+\frac{Af(r)}{\varphi_\alpha(r)}\right)^p 
	\le\varphi_\alpha^p(r) \left( 1+(\kappa+1)p \frac{Af(r)}{\varphi_\alpha(r)} \right) \\
	&=& \varphi_\alpha^p(r) + A(\kappa+1) p \varphi_\alpha^{p-1}(r) f(r)
	\qquad \mbox{for all } r>r_{A\eps}.
  \eas
  Since $\varphi_\alpha$ satisfies (\ref{eq:varphi}),
by (\ref{1.2}) 
  we have
  \bas
	& & (v_{0A\eps})_{rr} + \frac{n-1}{r} (v_{0A\eps})_r + v_{0A\eps}^p
	= (\varphi_{\alpha})_{rr} + Af_{rr} +
\frac{n-1}{r}(\varphi_{\alpha})_{r} + A \frac{n-1}{r} f_r
	+ v_{0A\eps}^p \\
	& &\qquad \le (\varphi_{\alpha})_{rr} + \frac{n-1}{r}(\varphi_{\alpha})_{r} + \varphi_\alpha^p 
          + A \Big\{ f_{rr} + \frac{n-1}{r} f_r + (\kappa+1) \varphi_\alpha^{p-1} f \Big\} 
	= 0
\eas
for all $r>r_{A\eps}$,
  which completes the proof of (\ref{2.3}).

  To conclude, it only remains to observe that if we let $r_A:=r_3(A)$ then 
  (\ref{2.10}) implies (\ref{2.01}), and that the claimed regularity
  properties of $v_{0A\eps}$ and (\ref{2.4}) are immediate consequences of the smoothness of $\varphi_\alpha$
  and the definition of $r_{A\eps}$.
\qed
We subsequently consider the radial version of (\ref{eq:main}), that is, we
investigate nonnegative solutions of
\be{0r}
        \left\{ \begin{array}{ll}
        (v^p)_t = v_{rr} + \frac{n-1}{r} v_r + v^p, \qquad & r>0, \ t>0, \\
        v(r,0)=v_0(r), \qquad & r\ge 0,
        \end{array} \right.
\ee
and to this end we introduce the operator $\parab$ defined by
  \be{parab}
        \parab v:=(v^p)_t - v_{rr} - \frac{n-1}{r} v_r - v^p.
  \ee
Then given radial initial data above $\varphi_\alpha$ but suitably close to $\varphi_\alpha$ asymptotically,
we may compare the corresponding solution $v$ with certain solutions of (\ref{0r}), emanating from appropriate initial data 
taken from Lemma \ref{lem2}, to show that the deviation $v-\varphi_\alpha$ essentially maintains its spatial decay
throughout evolution.
We can moreover make sure that $v$ approaches $\varphi_\alpha$ in the large time limit, yet without any information
on the rate of convergence.
\begin{lem}\label{lem3}
  Let $\alpha>0$, and assume that $v_0$ satisfies {\rm (\ref{sing})} and
  \be{3.1}
	\varphi_\alpha(r) \le v_0(r) \le \varphi_\alpha(r) + br^{-\gamma} 
	\qquad \mbox{for all } r>1
  \ee
  with some $b>0$ and $\gamma \in (\nu+\lambda_1,(n-2)/2)$.
  Then for the solution $v$ of {\rm (\ref{0r})} we have
  \be{3.2}
	\sup_{r\ge 0} \Big| v(r,t)-\varphi_\alpha(r)\Big| \to 0
	\qquad \mbox{as } t\to\infty.
  \ee
  Moreover, there exists $C>0$ such that
  \be{3.22}
	v(r,t) \le \varphi_\alpha(r) + C(r+1)^{-\gamma}
	\qquad \mbox{for all $r\ge 0$ and } t \ge 0.
  \ee
\end{lem}
\proof
  It can easily be checked that since $\gamma\in (\nu+\lambda_1,(n-2)/2)$, the number
$\kappa$ introduced in (\ref{kappa_gamma})
  is positive and satisfies (\ref{1.1}), and that with $\gamma(\kappa)$ as in (\ref{1.4}) we have $\gamma=\gamma(\kappa)$.
  Therefore, Lemma~\ref{lem1} applies to yield $c_1>0$ such that the solution $f=f_{\alpha,\kappa}$ of (\ref{1.2})
  satisfies
  \be{3.3}
	f(r)\ge c_1 r^{-\gamma}
	\qquad \mbox{for all } r>1.
  \ee
  Taking $b$ and $A_0$ from (\ref{3.1}) and Lemma~\ref{lem2}, respectively, we now fix $A>A_0$ such that
  \be{3.4}
	A\ge \frac{b}{c_1},
  \ee
  and let $r_A>1$ be as provided by Lemma~\ref{lem2}. 
  Then since $v_0(r) < Lr^{-\nu}$ for all $r>0$ by (\ref{sing}), the function $\chi$ given by
  \bas
	\chi(r):=\bigg( \frac{L}{v_0(r)} \bigg)^\frac{p-1}{2}-r, \qquad r\in [0,r_A],
  \eas
  is positive, and hence $\eps_0:=\min_{r\in [0,r_A]} \chi(r)$ satisfies $\eps_0>0$. 
  This enables us to finally choose $\eps\in (0,1)$ such that $\eps<\eps_0$, and let $r_{A\eps} \in (1,r_A)$ and
  $v_{0A\eps}$ be as given by Lemma~\ref{lem2}.

  Then since $r_{A\eps}>1$, (\ref{3.1}), (\ref{3.3}) amd (\ref{3.4}) imply that
  \bas
	v_0(r) - v_{0A\eps}(r)
	\le \Big(\varphi_\alpha(r)+ br^{-\gamma}\Big) - \Big(\varphi_\alpha(r)+ Af(r)\Big) 
	\le br^{-\gamma} - Ac_1 r^{-\gamma} 
	\le 0
\eas
	for all $r\ge r_{A\eps}$,
  whereas the inequality $r_{A\eps}<r_A$ in combination with our choice of $\chi$ and $\eps$ ensures that
  \bas
	v_0(r) = L\Big(r+\chi(r)\Big)^{-\nu}
	\le L(r+\eps_0)^{-\nu}
	\le L(r+\eps)^{-\nu}
	= v_{0A\eps}(r)
	\qquad \mbox{for all } r<r_{A\eps}.
  \eas
  This means that if we let $\ov$ denote the solution of
  \be{3.7}
	\left\{ \begin{array}{l}
	(\ov^p)_t=\ov_{rr}+ \frac{n-1}{r} \ov_r + \ov^p, \qquad r>0, \ t\in (0,T), \\[1mm]
	\ov(r,0)=v_{0A\eps}(r), \qquad r>0,
	\end{array} \right.
  \ee
  defined up to its maximal existence time $T\in (0,\infty]$,
  then
  \be{3.70}
	v_0(r) \le v_{0A\eps}(r)=\ov(r,0)
	\qquad \mbox{for all } r\ge 0.
  \ee
  But the properties (\ref{2.3}) and (\ref{2.4}) entail that $[0,\infty)^2 \ni (r,t) \mapsto v_{0A\eps}(r)$
  is a supersolution of (\ref{0r}) in the natural weak sense, which implies (\cite{quittner_souplet}) that the solution
  $\ov$ of (\ref{3.7}) satisfies
  \be{3.77}
	\ov_t(r,t) \le 0
	\qquad \mbox{for all $r\ge 0$ and } t\in (0,T).
  \ee
  Since clearly $\ov(r,t) \ge \varphi_\alpha(r)$ for all $r\ge 0$ and $t\in (0,T)$ by comparison, this entails that
  actually $T=\infty$, and that hence
  \be{3.8}
	\ov(r,t) \searrow v_\infty(r) 
	\qquad \mbox{as } t\to\infty
  \ee
  with some limit function $v_\infty$ fulfilling $\varphi_\alpha(r) \le v_\infty(r) \le v_{0A\eps}(r)$ for all $r\ge 0$.
  By a straightforward parabolic regularity argument (\cite{LSU}) based on the two-sided bound
  $\varphi_\alpha(r)\le \ov(r,t)\le v_{0A\eps}(r)$, $(r,t)\in [0,\infty)^2$, the convergence in (\ref{3.8}) can be seen
  to take place in $C^0([0,\infty)) \cap C^2_{loc}([0,\infty))$, which implies that $v_\infty$ is a stationary solution
  of (\ref{0r}). This means that with $\alpha':=v_\infty(0)\ge \varphi_\alpha(0)=\alpha$ we must have
  $v_\infty\equiv \varphi_{\alpha'}$, so that for the verification of (\ref{3.2}) it remains to be shown that
  $\alpha'\le \alpha$.

  Indeed, if we had $\alpha'>\alpha$ then by (\ref{a_alpha}) we would have $a_\alpha>a_{\alpha'}$ and hence could pick
  $a$ and $a'$ such that $a_{\alpha'} < a' < a < a_\alpha$. By (\ref{phi_exp}) we could thus find $r_1>0$ fulfilling
  \bas
	\varphi_{\alpha'}(r) \ge Lr^{-\nu} - a' r^{-\nu-\lambda_1}
	\qquad \mbox{for all } r>r_1
  \eas
  and
  \bas
	\varphi_{\alpha}(r) \le Lr^{-\nu} - a r^{-\nu-\lambda_1}
	\qquad \mbox{for all } r>r_1,
  \eas
  whereas Lemma~\ref{lem1} in view of the fact that $\gamma>\nu+\lambda_1$ would allow us to fix $r_2>0$
  satisfying
  \bas
	Af(r) < (a-a') r^{-\nu-\lambda_1}
	\qquad \mbox{for all } r>r_2.
  \eas
  With $r_3:=\max\{r_{A\eps},r_1,r_2\}$ we would thus arrive at the absurd conclusion that
  \bas
	Lr^{-\nu} - a' r^{-\nu-\lambda_1}
	&\le& \varphi_{\alpha'}(r) \le v_{0A\eps}(r) = \varphi_\alpha(r) + Af(r) \\
	&<& Lr^{-\nu} - a r^{-\nu-1} + (a-a') r^{-\nu-\lambda_1} \\
	&=& Lr^{-\nu} - a' r^{-\nu-\lambda_1}
	\qquad \mbox{for all } r>r_3,
  \eas
  from which we infer that in fact $\alpha'=\alpha$ and that thus (\ref{3.2}) is valid.
  Finally, (\ref{3.22}) immediately results from (\ref{3.77}), (\ref{3.70}) and the upper estimate for $f$ in (\ref{1.3}).
\qed
Another comparison argument shows that under the hypotheses of the last lemma, the solution will furthermore
eventually lie below any of the equilibria which are larger than the asymptotic profile.
\begin{lem}\label{lem5}
  Let $v_0$ be such that {\rm (\ref{sing})} holds, and such that {\rm (\ref{3.1})} is valid with some $\alpha>0$, $b>0$ and 
  $\gamma \in (\nu+\lambda_1,(n-2)/2)$. Then for any $\alpha'>\alpha$ one can find $t_0\ge 0$
  such that the solution $v$ of {\rm (\ref{0r})} satisfies
  \be{5.1}
	v(r,t)\le\varphi_{\alpha'}(r)
	\qquad \mbox{for all $r\ge 0$ and } t\ge t_0.
  \ee
\end{lem}
\proof
  In accordance with Lemma~\ref{lem3}, let us fix $c_1>0$ such that
  \be{5.2}
	v(r,t) \le \varphi_\alpha(r) + c_1(r+1)^{-\gamma}
	\qquad \mbox{for all $r\ge 0$ and } t\ge 0.
  \ee
  Since $\alpha'>\alpha$, by (\ref{a_alpha}) we can moreover pick positive numbers $a$ and $a'$ such that
  $a_{\alpha'}<a'<a<a_\alpha$, and thereupon use (\ref{phi_exp}) to find $r_1>0$ such that
  \bas
	\varphi_{\alpha'}(r) \ge Lr^{-\nu} - a' r^{-\nu-\lambda_1}
	\qquad \mbox{for all } r>r_1
  \eas
  as well as
  \bas
	\varphi_{\alpha}(r) \le Lr^{-\nu} - a r^{-\nu-\lambda_1}
	\qquad \mbox{for all } r>r_1.
  \eas
  This shows that if we abbreviate $c_2:=a-a'$ then
  \be{5.3}
	\varphi_{\alpha'}(r)-\varphi_\alpha(r) \ge c_2 r^{-\nu-\lambda_1}
	\qquad \mbox{for all } r > r_1.
  \ee
  Now thanks to the fact that $\gamma>\nu+\lambda_1$, choosing 
  \bas
	r_2:=\max\bigg\{r_1,\Big(\frac{c_1}{c_2}\Big)^\frac{1}{\gamma-\nu-\lambda_1} \bigg\}
  \eas
  we see that (\ref{5.2}) and (\ref{5.3}) imply the inequality
  \bea{5.33}
	v(r,t) &\le& \varphi_\alpha(r) + c_1 r^{-\gamma} 
	\le \varphi_{\alpha'}(r) - c_2 r^{-\nu-\lambda_1} + c_1 r^{-\gamma} \nn\\
	&=& \varphi_{\alpha'}(r) - c_2 r^{-\nu-\lambda_1}  
		\left( 1 - \frac{c_1}{c_2} r^{-(\gamma-\nu-\lambda_1)} \right) \nn\\
	&\le& \varphi_{\alpha'}(r) - c_2 r^{-\nu-\lambda_1} 
		\left( 1 - \frac{c_1}{c_2} r_2^{-(\gamma-\nu-\lambda_1)} \right) \nn\\
	&=&  \varphi_{\alpha'}(r) 
	\qquad \mbox{for all $r>r_2$ and } t>0.
  \eea
  Next, using that $\varphi_{\alpha'}(r)>\varphi_\alpha(r)$ for all $r\ge 0$ we can find $c_3>0$ fulfilling
  \be{5.4}
	\varphi_{\alpha'}(r)-\varphi_\alpha(r) \ge c_3
	\qquad \mbox{for all } r\in [0,r_2],
  \ee
  whereas the convergence statement in Lemma~\ref{lem3} provides $t_0\ge 0$ such that
  \be{5.5}
	v(r,t) \le \varphi_\alpha(r) + c_3
	\qquad \mbox{for all $r\ge 0$ and } t\ge t_0.
  \ee
  Combining (\ref{5.4}) and (\ref{5.5}) we thus obtain that
  \bas
	v(r,t) \le \varphi_{\alpha'}(r)
	\qquad \mbox{for all $r\in [0,r_2]$ and } t\ge t_0,
  \eas
  which together with (\ref{5.33}) establishes (\ref{5.1}).
\qed

In a third comparison procedure, we can finally
establish a quantitative upper estimate of the form asserted in Theorem \ref{theo100}.
For the first time we shall use here comparison functions which deviate from $\varphi_\alpha$ in a separated manner.
\begin{lem}\label{lem4}
  Assume that $v_0$ satisfies {\rm (\ref{sing})}, and that there exist $\alpha>0$, $b>0$ and 
  $\gamma \in (\nu+\lambda_1,(n-2)/2)$ such that {\rm (\ref{3.1})} holds.
  Then there exists $C>0$ such that the solution $v$ of {\rm (\ref{0r})} satisfies
  \be{4.1}
	|v(r,t)-\varphi_\alpha(r)| \le C e^{-\kappa(\gamma) t}
	\qquad \mbox{for all $r\ge 0$ and } t\ge 0
  \ee
with $\kappa(\gamma)$ as given by {\rm (\ref{kappa_gamma})}.
\end{lem}
\proof
  We fix any $\beta>\alpha+1$ and let $f:=f_{\beta,\kappa}$ denote the 
corresponding solution of (\ref{1.2}).
  Then according to (\ref{1.3}) we can find $c_1>0$ and $c_2>0$ such that
  \be{4.3}
	c_1 r^{-\gamma} \le f(r) \le c_2 r^{-\gamma}
	\qquad \mbox{for all } r>1,
  \ee
  whence in particular by (\ref{3.22}) and the positivity of $f$ we can pick $A>0$ such that
  \be{4.4}
	v(r,t) \le \varphi_\alpha(r) + Af(r)
	\qquad \mbox{for all $r\ge 0$ and } t\ge 0.
  \ee
  Next, since $\beta-1>\alpha$ we may apply Lemma~\ref{lem7} with $B:=0$ to find $r_1>1$ and $c_3>0$ satisfying
  \be{4.11}
	\varphi_{\beta-1}^{p-1}(r) - \Big( \varphi_\alpha(r)+Af(r) \Big)^{p-1}
	\ge c_3 r^{-2-\lambda_1}
	\qquad \mbox{for all } r>r_1,
  \ee
  whereas Lemma~\ref{lem9} provides $c_4>0$ such that
  \be{4.110}
	\varphi_\alpha^{p-2}(r) f(r) \le c_4 (r+1)^{\nu-2-\gamma}
	\qquad \mbox{for all } r\ge 0,
  \ee
  which clearly entails that
  \be{4.12}
	\varphi_\alpha^{p-2}(r) f(r) \le c_4 r^{\nu-2-\gamma}
	\qquad \mbox{for all } r>0,
  \ee
  because $\gamma+2>\nu$.

  Let us now pick $c_5>0$ and then $c_6>0$ such that
  \be{4.111}
	\varphi_\alpha(r) \ge c_5 r^{-\nu}
	\qquad \mbox{for all } r>1,
  \ee
  and such that with $z_1:=c_2 A/c_5$ we have
  \be{4.121}
	(1+z)^p \le 1+pz + c_6 z^2
	\qquad \mbox{for all } z\in [0,z_1].
  \ee
  Then since
  \be{4.122}
	\frac{Af(r)}{\varphi_\alpha(r)} \le \frac{c_2 A}{c_5} r^{-(\gamma-\nu)}
	\le \frac{c_2 A}{c_5}=z_1
	\qquad \mbox{for all } r>1
  \ee
  by (\ref{4.3}) and (\ref{4.111}), it follows from (\ref{4.11}) and (\ref{4.12}) that
  \bas
	\frac{\kappa p \Big\{\varphi_\beta^{p-1}(r) - \Big(\varphi_\alpha(r)+Af(r)\Big)^{p-1}\Big\}}
		{c_6 A \varphi_\alpha^{p-2}(r) f(r)}
	\ge \frac{\kappa p c_3}{c_4 c_6 A} r^{\gamma-\nu-\lambda_1}
	\qquad \mbox{for all } r>r_1.
  \eas
  As $\gamma>\nu+\lambda_1$, we can thus choose $r_2>r_1$ such that
  \be{4.13}
	\kappa p \Big\{\varphi_\beta^{p-1}(r) - \Big(\varphi_\alpha(r)+Af(r)\Big)^{p-1}\Big\} f(r)
	\ge c_6 A \varphi_\alpha^{p-2}(r) f^2(r)
	\qquad \mbox{for all } r>r_2.
  \ee
  With this value of $r_2$ fixed, by (\ref{phi_exp}) we can easily find $z_2>0$ such that
  \be{4.14}
	\frac{\varphi_\beta(r)}{\varphi_\alpha(r)} \le z_2
	\qquad \mbox{for all } r\in [0,r_2]
  \ee
  and thereupon let $c_7>0$ be large enough satisfying
  \be{4.15}
	(1+z)^p \le 1 + pz + c_7 z^2
	\qquad \mbox{for all }z\in [0,z_2].
  \ee
  Recalling (\ref{4.110}), we then take $c_8>0$ such that
  \be{4.16}
	c_7 \varphi_\alpha^{p-2}(r) f(r) \le c_8
	\qquad \mbox{for all } r\in [0,r_2],
  \ee
  and use that $\varphi_\beta>\varphi_{\beta-1}$ and $f>0$ on $[0,\infty)$ to obtain $c_9>0$ fulfilling
  \be{4.17}
	p \Big[\varphi_\beta^{p-1}(r)-\varphi_{\beta-1}^{p-1}(r) \Big] f(r) \ge c_9
	\qquad \mbox{for all } r\in [0,r_2].
  \ee
  Fixing $\delta\in (0,1)$ suitably small such that
  \be{4.18}
	c_8\delta \le c_9,
  \ee
  by an argument involving continuous dependence for the initial-value problem 
(\ref{eq:varphi})
  we are now able to find some $\alpha'\in (\alpha,\beta-1)$ sufficiently close to $\alpha$ such that
  \be{4.19}
	\varphi_{\alpha'}(r)-\varphi_\alpha(r) \le \delta
	\qquad \mbox{for all }r\in [0,r_2].
  \ee
  Finally, Lemma~\ref{lem5} allows us to choose $t_0>0$ large enough such that
  \be{4.20}
	v(r,t) \le \varphi_{\alpha'}(r)
	\qquad \mbox{for all $r\in [0,r_2]$ and } t\ge t_0.
  \ee
  We proceed to define a comparison function $\ov$ on $[0,\infty) \times [t_0,\infty)$ by letting
  \bas
	\ov(r,t):=\min \Big\{ \varphi_{\alpha'}(r) \, , \, \varphi_\alpha(r) + f(r)g(t) \Big\}
	\qquad \mbox{for $r\ge 0$ and } t\ge t_0,
  \eas
  where $g(t):=A e^{-\kappa(t-t_0)}$ for $t\ge t_0$.
  Then inside the set
  \bas
	Q:= \Big\{ (r,t) \in (0,\infty) \times (t_0,\infty) \ \Big| \ f(r)g(t) < \varphi_{\alpha'}(r) \Big\},
  \eas
  we compute
  \bas
	\parab \ov(r,t)
	&=& p\ov^{p-1} \ov_t - \ov_{rr} - \frac{n-1}{r} \ov_r - \ov^p \\[2mm]
	&=& - \kappa p  \Big[ \varphi_\alpha(r) - f(r)g(t) \Big]^{p-1} f(r)g(t) 
	 -\varphi_{\alpha rr}(r) - \frac{n-1}{r} \varphi_{\alpha r}(r) \\
	& &- g(t)  \Big[ f_{rr}(r)+\frac{n-1}{r}f_r(r) \Big] 
	 - \Big[ \varphi_\alpha(r) + f(r)g(t) \Big]^p \\[2mm]
	&=& g(t) \Bigg\{ - f_{rr}(r) - \frac{n-1}{r} f_r(r)
	- \kappa p \Big[ \varphi_\alpha(r) + f(r)g(t) \Big]^{p-1} f(r) \\
	& & \hspace*{10mm}
	- \frac{1}{g(t)} \varphi_\alpha^p(r) \bigg[ \Big(1+\frac{f(r)g(t)}{\varphi_\alpha(r)} \Big)^p-1\bigg] \Bigg\}
	\qquad \mbox{for all } (r,t)\in Q,
  \eas
  whence (\ref{1.2}) shows that
  \bea{4.200}
	\frac{1}{g(t)}  \parab \ov(r,t)
	&=& (\kappa+1)p \varphi_\beta^{p-1}(r) f(r) 
	-\kappa p \Big[ \varphi_\alpha(r) + f(r)g(t)\Big]^{p-1} f(r) \nn\\
	& & - \frac{1}{g(t)} \varphi_\alpha^p(r)  \bigg[ \Big(1+\frac{f(r)g(t)}{\varphi_\alpha(r)} \Big)^p-1\bigg]
	\qquad \mbox{for all } (r,t)\in Q.
  \eea
  Hence, if $(r,t)\in Q$ is such that $r>r_2$, then by (\ref{4.122}), (\ref{4.121}) and the fact that $g\le A$ we see that
  \bas
	\frac{1}{g(t)} \parab \ov(r,t)
	&\ge& (\kappa+1)p \varphi_\beta^{p-1}(r) f(r)
	-\kappa p \Big[ \varphi_\alpha(r)+Af(r)\Big]^{p-1} f(r) \\
	& & - \frac{1}{g(t)}\varphi_\alpha^p(r)  \Big\{ p \frac{f(r)g(t)}{\varphi_\alpha(r)}
		+ c_6  \frac{f^2(r) g^2(t)}{\varphi_\alpha^2(r)} \Big\} \\[2mm]
	&=& \kappa p \bigg\{ \varphi_\beta^{p-1}(r) - \Big[ \varphi_\alpha(r)+Af(r)\Big]^{p-1} \bigg\} f(r) \\
	& & + p \Big(\varphi_\beta^{p-1}(r)-\varphi_\alpha^{p-1}(r)\Big) f(r)
	- c_6 \varphi_\alpha^{p-2}(r) f^2(r) g(t).
  \eas
  Again since $g\le A$, and since $\varphi_\beta\ge \varphi_\alpha$, from (\ref{4.13}) we therefore obtain that
  \be{4.21}
	\parab \ov(r,t) \ge 0
	\qquad \mbox{whenever $(r,t)\in Q$ is such that } r>r_2.
  \ee
  On the other hand, if $(r,t)\in Q$ satisfies $r\le r_2$, then in particular
  \be{4.22}
	f(r)g(t) < \varphi_{\alpha'}(r) - \varphi_\alpha(r),
  \ee
  so that (\ref{4.14}) and the fact that $\alpha'<\beta$ imply that
  \bas
	\frac{f(r)g(t)}{\varphi_\alpha(r)} \le \frac{\varphi_\beta(r)-\varphi_\alpha(r)}{\varphi_\alpha(r)} \le z_2.
  \eas
  We may thus invoke (\ref{4.15}) and once more rely on (\ref{4.22}) to conclude from (\ref{4.200}) that for such $(r,t)$
  we have
  \bas
	\frac{1}{g(t)} \parab \ov(r,t)
	&\ge& (\kappa+1)p \varphi_\beta^{p-1}(r) f(r)
	-\kappa p  \Big[ \varphi_\alpha(r)+f(r)g(t)\Big]^{p-1}  f(r) \\
	& & - \frac{\varphi_\alpha^p(r)}{g(t)}  \bigg\{ p\frac{f(r)g(t)}{\varphi_\alpha(r)}
	- c_7  \frac{f^2(r)g^2(t)}{\varphi_\alpha^2(r)} \bigg\} \\[2mm]
	&=& \kappa p  \bigg\{ \varphi_\beta^{p-1}(r) - \Big[\varphi_\alpha(r)+f(r)g(t)\Big]^{p-1} \bigg\}  f(r) \\
	& & + p \Big(\varphi_\beta^{p-1}(r) -\varphi_\alpha^{p-1}(r)\Big)  f(r)
	- c_7 \varphi_\alpha^{p-2}(r) f^2(r)g(t) \\[2mm]
	&\ge& \kappa p  \Big\{ \varphi_\beta^{p-1}(r)-\varphi_{\alpha'}^{p-1}(r) \Big\}  f(r) \\
	& & + p \Big(\varphi_\beta^{p-1}(r) -\varphi_\alpha^{p-1}(r)\Big)  f(r)
	- c_7 \varphi_\alpha^{p-2}(r)  \Big(\varphi_{\alpha'}(r)-\varphi_\alpha(r)\Big)  f(r).
  \eas
  Since $\alpha'<\beta$ and $\alpha\le \beta-1$, in light of (\ref{4.16}), (\ref{4.17}), (\ref{4.19}), (\ref{4.18})
  and the fact that $\delta<1$, this yields
  \bas
	\frac{1}{g(t)}  \parab \ov(r,t)
	&\ge& p \Big(\varphi_\beta^{p-1}(r) -\varphi_{\beta-1}^{p-1}(r)\Big)  f(r)
	- c_7 \varphi_\alpha^{p-2}(r)  \delta f(r) \\[2mm]
	&\ge& c_9 - \delta c_8 
	\ge 0
	\qquad \mbox{for all $(r,t)\in Q$ with } r\le r_2.
  \eas
  In view of the equilibrium property of $\varphi_{\alpha'}$, this shows that $\ov$ is a supersolution of (\ref{0r})
  in $(0,\infty) \times (t_0,\infty)$.
  As (\ref{4.4}) and (\ref{4.20}) ensure that
  \bas
	v(r,t_0) \le \ov(r,t_0)
	\qquad \mbox{for all } r\ge 0,
  \eas
  by parabolic comparison we infer that $v(r,t)\le \ov(r,t)$ for all $r\ge 0$ and $t\ge t_0$, which evidently yields
  (\ref{4.1}) because of the boundedness of $f$.
\qed
\mysection{Convergence from below: upper bound for the rate}
For initial data below $\varphi_\alpha$, a quantitative convergence result can be derived by using an argument
which is based on a single comparison procedure, and which is thus somewhat simpler than the reasoning in the 
previous section.
\begin{lem}\label{lem10}
  Assume that $\alpha>0$, and that $v_0$ is nonnegative such that
  \be{10.1}
	\varphi_\alpha(r) \ge v_0(r) \ge \varphi_\alpha(r) - b (r+1)^{-\gamma}
	\qquad \mbox{for all } r\ge 0
  \ee
  with some $b>0$ and $\gamma\in (\nu+\lambda_1,(n-2)/2)$.
  Then there is $C>0$ such that the solution $v$ of {\rm (\ref{0r})} satisfies
  \be{10.2}
	|v(r,t)-\varphi_\alpha(r)| \le C e^{-\kappa(\gamma)t}
	\qquad \mbox{for all } t\ge 0,
  \ee
 with $\kappa(\gamma)$ as in {\rm (\ref{kappa_gamma})}.
\end{lem}
\proof
  As $\varphi_\alpha$ is an equilibrium of (\ref{0r}), the first inequality in (\ref{10.1}) along with
  a parabolic comparison shows that
  \bas
	v(r,t)\le \varphi_\alpha(r)
	\qquad \mbox{for all $r\ge 0$ and } t\ge 0,
  \eas
  so that we only need to establish a lower bound for $v$.
  To this end, we let $f=f_{\alpha,\kappa}$ be as given by Lemma~\ref{lem1}, so that from (\ref{1.3}) we know that
  \bas
	f(r)\ge c_1 r^{-\gamma} \ge c_1(r+1)^{-\gamma}
	\qquad \mbox{for all } r>1
  \eas
  with some $c_1>0$. Since $f$ is positive on $[0,1]$, we can thus pick $c_2>0$ such that
  \be{10.3}
	f(r)\ge c_2(r+1)^{-\gamma} 
	\qquad \mbox{for all } r\ge 0,
  \ee
  and let $A:=b/c_2$. Then again writing $g(t):=Ae^{-\kappa t}$, $t\ge 0$,
  with $\kappa:=\kappa(\gamma)$, we see that
  \bas
	\uv(r,t):=\max \Big\{0 \, , \, \varphi_\alpha(r)-f(r)g(t)\Big\},
	\qquad r\ge 0, \ t\ge 0,
  \eas
  satisfies 
  $\uv(r,0)\le v_0(r)$ for all $r\ge 0$.
  This is obvious whenever $\uv(r,0)=0$, while otherwise (\ref{10.3}) and (\ref{10.1}) assert that
  \bas
	\uv(r,0) =\varphi_\alpha(r)-Af(r)
	\le \varphi_\alpha(r) - Ac_2 (r+1)^{-\gamma}
	=\varphi_\alpha(r) - b(r+1)^{-\gamma}
	\le v_0(r).
  \eas
  In order to show that $\uv$ is a subsolution of (\ref{0r}) in $(0,\infty)^2$, we evidently only need to consider
  points $(r,t)$ where $\uv(r,t)>0$, at which we compute
  \bas
	\parab \uv(r,t)
	&=& \kappa p  \Big\{\varphi_\alpha(r)-f(r)g(t)\Big\}^{p-1}  f(r)g(t) 
	 - (\varphi_{\alpha})_{rr}(r) - \frac{n-1}{r} (\varphi_{\alpha})_{r}(r) \\
	& & - g(t) \Big\{ f_{rr}(r)+\frac{n-1}{r}f_r(r)\Big\} 
	 - \Big\{\varphi_\alpha(r)-f(r)g(t)\Big\}^p.
  \eas
  By monotonicity and convexity, respectively, we see that at any such point we have
  \bas
	\Big\{\varphi_\alpha(r)-f(r)g(t)\Big\}^{p-1} \le \varphi_\alpha^{p-1}(r)
  \eas
  and
  \bas
	\Big\{\varphi_\alpha(r)-f(r)g(t)\Big\}^p \ge \varphi_\alpha^p(r) + p\varphi_\alpha^{p-1}(r)  f(r)g(t),
  \eas
  so that since $\varphi_\alpha$ is a solution of
  (\ref{eq:varphi}), at all those $(r,t)$ we obtain
  \bas
	\parab \uv(r,t)
	\le g(t)  \Big\{\kappa p \varphi_\alpha^{p-1}(r) f(r) - f_{rr}(r)-\frac{n-1}{r}f_r(r)
	+ p\varphi_\alpha^{p-1}(r) f(r)\Big\} =0
  \eas
  according to (\ref{1.2}). 
  Using the comparison principle, we thus conclude that
  $\uv(r,t)\le v(r,t)$ for all $r\ge 0$ and $t\ge 0$, whence
  \bas
	\varphi_\alpha(r)-v(r,t) \le \varphi_\alpha(r) - \uv(r,t) \le f(r)g(t)
	\qquad \mbox{for all $r\ge 0$ and } t\ge 0.
  \eas
  Since $f$ is bounded, by definition of $g$ this yields the desired
 inequality and thereby proves (\ref{10.2}).
\qed

Now a straightforward combination of Lemma \ref{lem10} and Lemma \ref{lem4} provides the claimed upper estimate
on the convergence rate.\abs
{\sc Proof of Theorem~\ref{theo100}.} \quad
If $v_0$ is as in Theorem~\ref{theo100} then one can find functions
$\uv_0,\ov_0$,
\[
\uv_0(|x|)\le\min\{v_0(x),\varphi_\alpha(|x|)\}, \qquad x\in\R^n,
\]
and
\[
\max\{v_0(x),\varphi_\alpha(|x|)\}\le \ov_0(|x|)\le L|x|^{-\nu}, \qquad
x\in\R^n,\quad x\not=0,
\]
such that $\uv_0,\ov_0$ satisfy the assumptions of Lemma~\ref{lem10},
Lemma~\ref{lem4}, respectively. By comparison, these two lemmata yield then
the claim.
\qed
\mysection{Convergence from above: lower bound for the rate}
In this section we shall prove the first statement in Theorem \ref{theo200}. 
For this purpose, we employ comparison functions similar to those used in Lemma \ref{lem4} and Lemma \ref{lem10}
to establish the following.
\begin{lem}\label{lem6}
  Let $\alpha>0$, and assume that $v_0$ satisfies {\rm (\ref{sing})} as well
as
  \be{6.1}
        v_0(r) \ge \varphi_\alpha(r) + b (r+1)^{-\gamma}
        \qquad \mbox{for all } r\ge 0
  \ee
  with some $b>0$ and $\gamma\in (\nu+\lambda_1,(n-2)/2)$.
  Then there exists $c>0$ such that the solution $v$ of {\rm (\ref{0r})}
  satisfies
  \be{6.2}
        v(0,t) - \varphi_\alpha(0) \ge c e^{-\kappa(\gamma)t}
        \qquad \mbox{for all } t\ge 0,
  \ee
 where $\kappa(\gamma)$ is as in {\rm (\ref{kappa_gamma})}.
\end{lem}
\proof
  We let $f:=f_{\alpha,\kappa}$ be as given by Lemma~\ref{lem1}, and then
obtain from (\ref{1.3}), (\ref{6.1}) and the
  boundedness of $f$ that there exists $A>0$ such that
  \bas
        v_0(r) \ge \varphi_\alpha(r) + Af(r)
        \qquad \mbox{for all } r\ge 0.
  \eas
  This means that $v$ initially dominates the function $\uv$ defined on
  $[0,\infty)^2$ by setting
  \bas
        \uv(r,t):=\varphi_\alpha(r) + f(r)g(t)
        \quad \mbox{for $r\ge 0$ and $t\ge 0$, \quad with \quad}
        g(t):=Ae^{-\kappa t}
        \quad \mbox{for } t\ge 0
  \eas
  and $\kappa:=\kappa(\gamma)$. Then
  \bea{6.3}
        \parab \uv(r,t)
        &=& - \kappa p  \Big\{\varphi_\alpha(r)+f(r)g(t)\Big\}^{p-1}
f(r)g(t)
         - (\varphi_{\alpha})_{rr}(r)-\frac{n-1}{r}(\varphi_{\alpha})_{r}(r)
         - \nn\\
         & &- g(t)  \Big\{f_{rr}(r)+\frac{n-1}{r} f_r(r)\Big\} 
        - \Big\{\varphi_\alpha(r)+f(r)g(t)\Big\}^p
\eea
for all $r>0$ and $t>0$,
  where by monotonicity and convexity of $[0,\infty)\ni z \mapsto z^{p-1}$ 
we can estimate
  \bas
        \Big\{\varphi_\alpha(r)+f(r)g(t)\Big\}^{p-1} \ge
\varphi_\alpha^{p-1}(r)
  \eas
  and 
  \bas
        \Big\{\varphi_\alpha(r)+f(r)g(t)\Big\}^p
        \ge \varphi_\alpha^p(r) + p\varphi_\alpha^{p-1}(r)  f(r)g(t)
  \eas
  for all $r>0$ and $t>0$. Since $\varphi_\alpha$ satisfies
(\ref{eq:varphi}),
  (\ref{6.3}) implies that
  \bas  
        \parab \uv (r,t)
        \le g(t)  \Big\{ -\kappa p\varphi_\alpha^{p-1}(r)f(r)
        - f_{rr}(r)-\frac{n-1}{r}f_r(r) - p\varphi_\alpha^{p-1}(r)f(r)\Big\}
        = 0
\eas
for all $r>0$ and $t>0$
  due to (\ref{1.2}).  
  The comparison principle thus shows that $\uv(r,t) \le v(r,t)$ for all
$r\ge 0$ and $t\ge 0$, whence in particular
  \bas
        v(0,t)\ge \uv(0,t) =\varphi_\alpha(0)+Ae^{-\kappa t}
        \qquad \mbox{for all }  t\ge 0,
  \eas
  because $f(0)=1$.
\qed

From this we immediately obtain the first conclusion of Theorem \ref{theo200}.\abs
{\sc Proof of Theorem~\ref{theo200}}~(i). \quad
If $v_0$ is as in Theorem~\ref{theo200}~(i) then one can find a function
$\uv_0\le v_0$ satisfying the assumptions of Lemma~\ref{lem6}.
By comparison, Lemma~\ref{lem6} yields then
the claim.
\qed
\mysection{Convergence from below: lower bound for the rate}
We shall next verify that the convergence statement in Theorem \ref{theo100} indeed yields the precise convergence rate
also for solutions approaching their limit from below. 
Paralleling Lemma \ref{lem6}, the following lemma provides a technical preparation for this.
\begin{lem}\label{lem8}
  Let $\alpha>0$ and $v_0$ be nonnegative and such that
  \be{8.1}
	v_0(r) \le \varphi_\alpha(r) - b (r+1)^{-\gamma}
	\qquad \mbox{for all } r\ge 0
  \ee
  with some $b>0$ and $\gamma\in (\nu+\lambda_1,(n-2)/2)$.
  Then there exists $c>0$ such that for the solution $v$ of {\rm (\ref{0r})} we have
  \be{8.2}
	\varphi_\alpha(0) - v(0,t) \ge c e^{-\kappa(\gamma)t}
	\qquad \mbox{for all } t\ge 0,
  \ee
 where $\kappa(\gamma)$ is as in {\rm (\ref{kappa_gamma})}.
\end{lem}
\proof
  Let us pick $\beta\in (0,\alpha)$ and take $f:=f_{\beta,\kappa}$ from Lemma~\ref{lem1}. 
  Then Lemma~\ref{lem7} yields $r_1\ge 1$ and $c_1>0$ such that
  \be{8.3}
	\Big\{\varphi_\alpha(r)-f(r)\Big\}^{p-1} - \varphi_\beta^{p-1}(r) \ge c_1 r^{-2-\lambda_1}
	\qquad \mbox{for all } r>r_1,
  \ee
  and from Lemma~\ref{lem9} we obtain $c_2>0$ such that
  \be{8.4}
	\varphi_\alpha^{p-2}(r)  f(r) \le c_2 (r+1)^{\nu-2-\gamma}
	\qquad \mbox{for all } r>0.
  \ee
  Moreover, since $f$ is bounded and $\varphi_\alpha$ is positive with
  \bas
	\frac{f(r)}{\varphi_\alpha(r)} \le c_3 r^{-\gamma+\nu}
	\qquad \mbox{for all } r>1
  \eas
  and some $c_3>0$ by (\ref{1.3}) and (\ref{phi_exp}), due to the fact that $\gamma>\nu$ we can fix $c_4>0$
  satisfying
  \be{8.5}
	\frac{f(r)}{\varphi_\alpha(r)} \le c_4
	\qquad \mbox{for all } r\ge 0.
  \ee
  Similarly, using (\ref{1.3}) and the boundedness of $f$ we can find $c_5>0$ such that
  \be{8.55}
	(r+1)^\gamma f(r)\le c_5
	\qquad \mbox{for all } r\ge 0.
  \ee
  We now pick $c_6>0$ fulfilling
  \be{8.6}
	(1-z)^p \le 1-pz+c_6 z^2
	\qquad \mbox{for all } z\in [0,1/2],
  \ee
  and finally choose $A\in (0,1]$ such that
  \be{8.8}
	A \le \min\left\{ \frac{1}{2c_4},\frac{b}{c_5},\frac{\kappa p c_1}{c_2 c_6} r_1^{\gamma-\nu-\lambda_1}\right\},
  \ee
  as well as
  \be{8.9}
	c_7:=\inf_{r\in (0,r_1)} \bigg\{ \Big\{\varphi_\alpha(r)-Af(r)\Big\}^{p-1} - \varphi_\beta^{p-1}(r) \bigg\} \, >0.
  \ee
 With these choices, we again define $g(t):=Ae^{-\kappa t}$ for $t\ge 0$
  where $\kappa:=\kappa(\gamma)$. Then we set
  \bas
	\ov(r,t):=\varphi_\alpha(r)-f(r)g(t)
	\qquad \mbox{for $r\ge 0$ and } t\ge 0,
  \eas
  and first observe that as a consequence of (\ref{8.5}) and (\ref{8.8}) we have
  \be{8.11}
	\frac{f(r)g(t)}{\varphi_\alpha(r)}
	\le \frac{Af(r)}{\varphi_\alpha(r)} \le c_4 A \le \frac{1}{2}
	\qquad \mbox{for all $r> 0$ and $t> 0$,}
  \ee
  whence in particular
  \bas
	\ov(r,t) \ge \frac{1}{2}\varphi_\alpha(r)
	\qquad \mbox{for all $r> 0$ and $t> 0$.}
  \eas
  Furthermore, (\ref{8.11}) allows us to apply (\ref{8.11}) in estimating
  \bas
	\ov^p(r,t)
	= \Big\{ \varphi_\alpha(r)-f(r)g(t)\Big\}^p
	\le \varphi_\alpha^p(r) - p\varphi_\alpha^{p-1}(r)  f(r)g(t) + c_6 \varphi_\alpha^{p-2}(r)  f^2(r)g^2(t),
  \eas
  so that using the equilibrium property of $\varphi_\alpha$ and (\ref{1.2}) we find that
  \bea{8.12}
	\parab \ov(r,t)
	&=& \kappa p  \Big\{ \varphi_\alpha(r)-f(r)g(t)\Big\}^{p-1}  f(r)g(t) 
	 - (\varphi_{\alpha})_{rr}(r)-\frac{n-1}{r}(\varphi_{\alpha})_{r}(r) \nn\\
	& &- g(t) \Big\{ f_{rr}(r)+\frac{n-1}{r}f_r(r)\Big\} 
	 - \ov^p(r,t) \nn\\
	&\ge& \kappa p  \Big\{ \varphi_\alpha(r)-f(r)g(t)\Big\}^{p-1}  f(r)g(t) 
	 - g(t) \Big\{ f_{rr}(r)+\frac{n-1}{r}f_r(r)\Big\} \nn\\
	& & + p\varphi_\alpha^{p-1}(r)  f(r)g(t) - c_6 \varphi_\alpha^{p-2}(r)  f^2(r)g^2(t) \nn\\
	&=& f(r)g(t)  \Big\{
	\kappa p  \left[ \left\{ \varphi_\alpha(r) - f(r)g(t)\right\}^{p-1} - \varphi_\beta^{p-1}(r)\right] \nn\\
	& & 
	+ p \left[\varphi_\alpha^{p-1}(r)-\varphi_\beta^{p-1}(r) \right] 
	- c_6 \varphi_\alpha^{p-2}(r) f(r)g(t) \Big\}
\eea
for all $r>0$ and $t>0$.
  Here
  \bas
	\varphi_\alpha^{p-1}(r)-\varphi_\beta^{p-1}(r)>0
	\qquad \mbox{for all } r>0,
  \eas
  because $\alpha>\beta$. Moreover, since $g(t) \le A$ for all $t\ge 0$, we see that (\ref{8.4}) implies
  \be{8.13}
	c_6 \varphi_\alpha^{p-2}(r) f(r)g(t)
	\le c_6 A \varphi_\alpha^{p-2}(r)f(r) 
	\le c_2 c_6 A (r+1)^{-2+\nu-\gamma}
	\qquad \mbox{for all } r>0,
  \ee
  and that
  \be{8.14}
	\kappa p  \bigg[ \Big\{ \varphi_\alpha(r) - f(r)g(t)\Big\}^{p-1} - \varphi_\beta^{p-1}(r)\bigg]
	\ge \kappa p  \bigg[ \Big\{\varphi_\alpha(r)-Af(r)\Big\}^{p-1}-\varphi_\beta^{p-1}(r)\bigg]
  \ee
  for all $r>0$.
  In particular, using (\ref{8.3}) and (\ref{8.8}) along with the fact that $A\le 1$, we obtain from
  (\ref{8.13}) and (\ref{8.14}) that
  \bea{8.15}
	\frac{\kappa p  \bigg[ \Big\{ \varphi_\alpha(r) - f(r)g(t)\Big\}^{p-1} - \varphi_\beta^{p-1}(r)\bigg]}
		{c_6 \varphi_\alpha^{p-2}(r) f(r)g(t)}
	\ge \frac{\kappa p  c_1 r^{-2-\lambda_1}}{c_2c_6 A (r+1)^{\nu-2-\gamma}} \ge 1
\eea
for all $r\ge r_1$ and $t>0,$
  because $\gamma>\nu+\lambda_1$. 
  On the other hand, for small $r$ we apply (\ref{8.9}) and (\ref{8.8}) to find that
  \bea{8.16}
	\frac{\kappa p  \bigg[ \Big\{ \varphi_\alpha(r) - f(r)g(t)\Big\}^{p-1} - \varphi_\beta^{p-1}(r)\bigg]}
		{c_6 \varphi_\alpha^{p-2}(r) f(r)g(t)}
	\ge \frac{\kappa p c_7}{c_2 c_6 A(r+1)^{\nu-2-\gamma}} \ge 1
\eea
for all $r\in (0,r_1)$ and $t>0$.
  From (\ref{8.12}), (\ref{8.15}) and (\ref{8.16}) we thus infer that
  \bas
	\parab \ov(r,t)\ge 0
	\qquad \mbox{for all $r>0$ and } t>0,
  \eas
  so that since for any $r\ge 0$ we have
  \bas
	\ov(r,0)
	= \varphi_\alpha(r)-Af(r)
	\ge \varphi_\alpha(r)-c_5 A (r+1)^{-\gamma}
	\ge \varphi_\alpha(r)-b (r+1)^{-\gamma}
	\ge v_0(r)
  \eas
  according to (\ref{8.55}), (\ref{8.8}) and (\ref{8.1}), by comparison we conclude that $\ov(r,t)\ge v(r,t)$ for
  all $r\ge 0$ and $t\ge 0$. In particular, this implies that
  \bas
	\varphi_\alpha(0)-v(0,t) \ge \varphi_\alpha(0)-\ov(0,t)
	=Ae^{-\kappa t}
	\qquad \mbox{for all } t\ge 0,
  \eas
  which proves (\ref{8.2}).
\qed

We can thereby complete the proof of Theorem \ref{theo200}.\abs
{\sc Proof of Theorem~\ref{theo200}}~(ii). \quad
If $v_0$ is as in Theorem~\ref{theo200}~(ii) then one can find a function
$\ov_0\ge v_0$ satisfying the assumptions of Lemma~\ref{lem8}.
By comparison, the claim follows from Lemma~\ref{lem8}.
\qed
\mysection{Instability of the steady states when $p< p_c$}
{\sc Proof of Proposition~\ref{prop:unst}}.\quad
The proof is based on intersection properties of the steady states similarly
as the proof of Theorem~1.14 in \cite{Gui-N-W:stab}. We recall that for $p_S\le p< p_c$
any two steady states intersect, see \cite{Wang}.

To prove (i) we may assume without loss of generality that
$v_0(x)>\varphi_\alpha(|x|)$ for $x\in\R^n$. Then we choose $\varepsilon >0$
small enough such that 
\[
v_*(|x|):=\max\{\varphi_\alpha(|x|),
\varphi_{\alpha+\varepsilon}(|x|)\}\le v_0(x),\qquad x\in\R^n.
\]
Then the solution $\uv$ of (\ref{eq:main}) with the initial condition
$\uv(x,0)=v_*(|x|)$ satisfies $(\uv)_t\ge 0$ in $\R^n\times(0,t_{max})$
where $t_{max}\in(0,\infty]$ is the maximal existence time. In fact,
$t_{max}=\infty$ for every solution of (\ref{eq:main}) with a bounded
initial function $v_0$ because $\Vert v_0\Vert_{L^\infty(\R^n)}e^{t/p}$ is a
supersolution.

Suppose the increasing function $\uv(0,t)=\Vert v(\cdot,t)\Vert_{L^\infty(\R^n)}$
has a finite limit as $t\to\infty$. Then $\uv(\cdot,t)$ converges to a
steady state that is bigger than $\varphi_\alpha$ which is a contradiction.

To prove (ii) we assume that
$v_0(x)<\varphi_\alpha(|x|)$ for $x\in\R^n$, and choose $\varepsilon >0$
small enough such that
\[
v^*(|x|):=\min\{\varphi_\alpha(|x|),
\varphi_{\alpha-\varepsilon}(|x|)\}\ge v_0(x),\qquad x\in\R^n.
\]
Then the solution $\ov$ of (\ref{eq:main}) with the initial condition
$\ov(x,0)=v^*(|x|)$ satisfies $(\uv)_t\le 0$ in $\R^n\times(0,T)$
where $T\in(0,\infty]$ is the maximal time such that $\ov>0$
in $\R^n\times(0,T)$. 

If the decreasing function $\ov(0,t)=\Vert
v(\cdot,t)\Vert_{L^\infty(\R^n)}$
had a positive limit as $t\to T$ then $\ov(\cdot,t)$ would converge to a
positive steady state smaller than $\varphi_\alpha$ which is a contradiction.
\qed
\noindent {\bf Acknowledgments.}
The first author was supported in part by the Slovak Research and
Development Agency under the contract
No. APVV-0134-10 and by the VEGA grant 1/0711/12.


\begin{thebibliography}{99}
%
\bibitem{AK}
\newblock G. Akagi and R. Kajikiya,
\newblock \emph{Stability analysis of asymptotic profiles for sign-changing
solutions to fast diffusion equations},
\newblock Manuscripta Math., {\bf 141} (2013), 559--587.

\bibitem{BH} 
\newblock J. G. Berryman and C. J. Holland,
\newblock \emph{Stability of the separable solution for fast diffusion},
\newblock Arch. Rat. Mech. Anal., \textbf{74} (1980), 379--388.

\bibitem{BBDGV} 
\newblock A. Blanchet, M. Bonforte, J. Dolbeault, G. Grillo and J. L.
V\'azquez,
\newblock \emph{Asymptotics of the fast diffusion equation via entropy
estimates},
\newblock Arch. Rat. Mech. Anal., \textbf{191} (2009), 347--385.

\bibitem{BDGV} 
\newblock M. Bonforte, J. Dolbeault, G. Grillo and J. L. V\'azquez,
\newblock \emph{Sharp rates of decay of solutions to the nonlinear fast
diffusion equation via functional inequalities},
\newblock Proc. Nat. Acad. Sciences, {\bf 107}
(2010), 16459--16464.

\bibitem{BGV} 
\newblock M. Bonforte, G. Grillo and J. L. V\'azquez,
\newblock \emph{Special fast diffusion with slow asymptotics. Entropy method
and flow on a Riemannian manifold},
\newblock Arch. Rat. Mech. Anal., \textbf{196} (2010), 631--680.

\bibitem{BGVb} 
\newblock M. Bonforte, G. Grillo and J. L. V\'azquez,
\newblock \emph{Behaviour near extinction for the Fast Diffusion Equation
on bounded domains},
\newblock J. Math. Pures Appl., {\bf 97} (2012), 1--38.

\bibitem{dPS} 
\newblock M. del Pino and M. S\'aez,
\emph{On the extinction profile for solutions of $u_t=\Delta u^{(N-2)/(N+2)}$},
\newblock Indiana Univ. Math. J., {\bf 50} (2001), 611--628.

\bibitem{FKW}
\newblock M. Fila, J. R. King and M. Winkler,
\newblock \emph{Rate of convergence to Barenblatt profiles for the
fast diffusion equation with a critical exponent}, 
\newblock J. London Math. Soc., to appear.

\bibitem{FVW} 
\newblock M. Fila, J. L. V\'azquez and M. Winkler,
\newblock \emph{A continuum of extinction
rates for the fast diffusion equation},  
\newblock Comm.~Pure Appl.~Anal., {\bf 10} (2011), 1129--1147.

\bibitem{FVWY}
\newblock M. Fila, J. L. V\'azquez, M. Winkler and E. Yanagida,
\emph{Rate of convergence to Barenblatt profiles for the fast diffusion
  equation},
\newblock Arch. Rat. Mech. Anal., {\bf 204} (2012), 599--625.

\bibitem{FWY}
\newblock M. Fila, M. Winkler and E. Yanagida,
\newblock \emph{Convergence rate for a parabolic equation with
supercritical nonlinearity},
\newblock J. Dynam. Differential Equations, {\bf 17} (2005),
 249--269.

\bibitem{GP}
\newblock V. A. Galaktionov and L. A. Peletier,
\newblock \emph{Asymptotic behaviour near finite-time extinction for the
fast diffusion equation},
\newblock Arch. Rat. Mech. Anal., {\bf 139} (1997), 83--98.

\bibitem{Gui-N-W:stab}
C.~Gui, W.-M. Ni, and X.~Wang,
\newblock \emph{On the stability and instability of positive steady states of a
  semilinear heat equation in {${\R}^n$}},
\newblock Comm. Pure Appl. Math., {\bf 45} (1992),
1153--1181.

\bibitem{HY}
\newblock M. Hoshino and Y. Yanagida,
\newblock \emph{Sharp estimates of the convergence rate for a semilinear 
parabolic equation with supercritical nonlinearity}, 
\newblock Nonlin. Anal. TMA, {\bf 69} (2008), 3136--3152.

\bibitem{K1}
\newblock J. R. King,
\newblock \emph{Self-similar behaviour for the equation of fast
nonlinear diffusion},
\newblock Phil. Trans. Roy. Soc. Lond. A, {\bf 343} (1993), 337--375.

\bibitem{K2}
\newblock J. R. King,
\newblock \emph{Asymptotic analysis of extinction behaviour in fast
nonlinear diffusion},
\newblock J. Eng. Math., {\bf 66} (2010), 65--86.

\bibitem{Kw}
\newblock Y. C. Kwong,
\newblock \emph{Asymptotic behavior of a plasma type equation with finite
extinction},
\newblock Arch. Rat. Mech. Anal., {\bf 104} (1988), 277--294.

\bibitem{LSU} 
\newblock O. A. Ladyzhenskaya, V. A. Solonnikov and N. N. Uraltseva,
\newblock ``Linear and Quasilinear Equations of Parabolic Type'',
\newblock Amer. Math. Soc., Providence, RI, 1968.

\bibitem{PZ} 
\newblock M. A. Peletier and H. Zhang,
\newblock \emph{Self-similar solutions of a fast diffusion equation that do
not conserve mass},
\newblock Diff. Int. Equations, {\bf 8} (1995), 2045--2064.

\bibitem{quittner_souplet} 
\newblock P. Quittner and Ph. Souplet,
\newblock ``Superlinear Parabolic Problems. Blow-up,
Global Existence and Steady States'',
\newblock Birkh\"auser Advanced Texts, Birkh\"auser, Basel, 2007.

\bibitem{SV}
\newblock G. Savar\'e and V. Vespri,
\newblock \emph{The asymptotic profile of solutions of a class of doubly
nonlinear equations},
\newblock Nonlin. Anal. TMA, {\bf 22} (1994), 1553--1565.

\bibitem{Vsmooth}
\newblock J. L. V\'azquez,
\newblock ``Smoothing and Decay Estimates for Nonlinear Diffusion
Equations'',
\newblock Oxford Lecture Notes in Maths. and its
Applications, vol.~33, Oxford University Press, Oxford, 2006.

\bibitem{Wang}
X.~Wang,
\newblock \emph{On the Cauchy problem for reaction-diffusion equations},
\newblock Trans. Amer. Math. Soc., {\bf 337} (1993), 549--590.

%
\end{thebibliography}
\end{document}